\documentclass{amsart}
\usepackage{graphicx}
\newtheorem{thm}{Theorem}[section]

\newtheorem{lem}[thm]{Lemma}

\theoremstyle{definition}

\theoremstyle{remark}
\newtheorem{rem}[thm]{Remark}
\numberwithin{equation}{section}

\begin{document}

\title[Statistical analysis of loss queueing systems]{Statistical analysis of single-server loss queueing systems}%
\author{Vyacheslav M. Abramov}%
\address{School of Mathematical Sciences, Monash University, Clayton Campus, Building 28, Level 4,
Wellington rd, VIC-3800, Australia}%
\email{vyacheslav.abramov@sci.monash.edu.au}%

\thanks{The author acknowledge with thanks the support of the Australian Research Council.}%

\begin{abstract}
In this article statistical bounds for certain output characteristics of the $M/GI/1/n$
and $GI/M/1/n$ loss queueing systems are derived on the basis of large samples of an input
characteristic of these systems, such as service time in the $M/GI/1/n$ queueing system or
interarrival time in the $GI/M/1/n$ queueing system. The analysis of this article is based
on application of Kolmogorov's statistics for empirical probability distribution functions.
\end{abstract}
\maketitle
\section{Introduction}
In theoretical problems of queueing theory, the input
characteristics such as interarrival and service time distributions
are assumed to be known. For example, if we speak about an $M/GI/1$
queueing system, we assume that the arrival process is Poisson with
rate $\lambda$, and service times are independent and identically
distributed random variables with a given probability distribution
function $B(x)$. In practice we, however, have only input data
characterizing arrival and departure processes, and our conclusion
about output characteristics depend on accuracy of approximation of
the aforementioned input characteristics.

In certain queueing systems some output characteristics can be
insensitive to the type of probability distribution function of a
service time. For example, in $M/GI/m/0$ queueing systems the
stationary state probabilities are defined by the Erlang-Sevastyanov
formulae, which are independent of the type of the probability
distribution function of a service time. In the $M/G/1$ queueing
system, in which arrival process is Poisson with rate $\lambda$, and
the expected service time is $b<\frac{1}{\lambda}$, the expected
length of a busy period is $\frac{b}{1-\lambda b}$, i.e. insensitive
to the type of probability distribution $B(x)$. Therefore, the
problem of estimation of the expected busy period and certain other
output characteristics, such as, for example, the expected number of
served customer during a busy period, reduces to those estimations
of the parameters $\lambda$ and $b$ of the queueing system.

The queueing systems, for which the output characteristics are
insensitive to the types of probability distributions of input
characteristics, are rather exceptions. In most cases output
characteristics depend on probability distribution functions, and
this dependence can be very complicated. So, the problem of
estimating the output characteristics of queueing systems is
generally difficult problem.

In the present paper we demonstrate the methods of statistical
analysis of certain output characteristics of single-server loss
queueing systems such as $M/GI/1/n$ and $GI/M/1/n$, which are
\textit{not} insensitive, but expressed specifically via transforms
of probability distribution function of interarrival or service
times. We estimate the output characteristics of these systems, such
as the expected busy period, expected numbers of lost and served
customers during a busy period for the $M/GI/1/n$ and the stationary
loss probability for the $GI/M/1/n$ queueing systems on the basis of
known information about input characteristics of these systems. The
concrete problem formulations are given later.

Statistics of queueing systems is a distinguished area of queueing
theory. The first publications in this area appeared long time ago
(see the textbook of Ivchenko, Kashtanov and Kovalenko
\cite{IvKashKov}, review papers of Bhat and Rao \cite{Bhat Subba
Rao} and Daley \cite{Daley 1976} and the references in these
sources). In the textbook \cite{IvKashKov} various traditional
methods of statistical inference to queueing problems have been
demonstrated. In \cite{Bhat Subba Rao}, a review of different
aspects of queueing systems, including identification of models,
parameters estimation by the maximum likelihood method and the
method of moments as well as estimates of mean value processes and
auto-covariance functions, hypothesis testing and other topics of
statistical analysis up to the publication date is made. In
\cite{Daley 1976}, a review of various aspects, including
statistical, concerning output or departure processes of $G/G/s/N$
queueing systems is made. For some recent publications in the area
of statistical analysis of queueing systems see also \cite{Acharya},
\cite{Basawa Bhat Lund}, \cite{Hansen Pitts}, \cite{Mandelbaum
Zeltin}, \cite{Novak Watson}. Nevertheless, despite their
importance, the papers on statistical inference of queueing systems
appears much rarely compared to many of those that use the methods
of stochastic analysis, optimization, control and asymptotic methods
of Mathematical analysis. Moreover, the methods of statistical
analysis of queueing systems known from the literature are
traditional.

The statistical analysis of the present paper is based on
application of Kolmogorov's statistics characterizing empirical
probability distribution function of an interarrival or service time
on the basis of large number of observations of these input
characteristics.

To our knowledge, Kolmogorov's statistics are never used in
statistical analysis of output characteristics of queueing systems.
In recent paper \cite{A1}, Kolmogorov's metric (which is associated
with one of Kolmogorov's statistics) is used for establishing the
bounds for the loss probability in certain queueing systems with
large buffers. In other papers \cite{A2} and \cite{A3} Kolmogorov's
metric is used for establishing conditions for the continuity in the
$M/M/1/n$ queueing system and, respectively, for the continuity of
non-stationary state probabilities in large closed queueing networks
with bottlenecks.

In the present paper we solve the following problems. Consider,
for instance, the $M/GI/1/n$ queueing system in which the arrival rate
$\lambda$ and the expected service time $b$ are assumed to be known, however,
the probability distribution function $B(x)$ of a service time is unknown.
On the basis of $N$ observations of service times we build an empirical probability
distribution function $B_{\mathrm{emp}}(x,N)$.
Let us denote
$$
\delta_N=\sup_{x>0}\big|B(x)-B_{\mathrm{emp}}(x,N)\big|.
$$
Using this statistic, we estimate the expected busy period, and
other significant characteristics such as the expected numbers of
served and lost customers during that busy period. We find
confidence intervals (ranges) for these characteristics based on the
confidence probability $P$. The similar ranges will be obtained with
the aid of other Kolmogorov's statistics mentioned later.

Recall that according to Kolmogorov's theorem
(see \cite{Kolmogorov} or \cite{Takacs 1967}, p. 170)
we have:
\begin{equation}\label{2}
\begin{aligned}
\lim_{N\to\infty}\mathrm{Pr}\left\{\delta_N<\frac{z}{\sqrt{N}}\right\}
&=K(z)\\ &=\begin{cases}\sum_{j=-\infty}^{+\infty}(-1)^j\mathrm{e}^{-2j^{2}z^{2}}, &\mbox{for} \ z>0,\\
0, &\mbox{for} \ z\leq0.
\end{cases}
\end{aligned}
\end{equation}
So,
for the range
\begin{equation*}\label{1}
\delta_N<\epsilon,
\end{equation*}
based on a chosen confidence probability $P$,
the value $\epsilon=\frac{z}{\sqrt{N}}$ depending on the large parameter $N$ can be chosen such that the equation $K(z)=P$ is satisfied.

Along with $\delta_N$, the other relevant Kolmogorov's statistics are as follows:
\begin{equation}\label{1'}
\delta_N^-=\sup_{x>0}\big[B(x)-B_{\mathrm{emp}}(x,N)\big],
\end{equation}
and
\begin{equation}\label{1''}
\delta_N^+=\sup_{x>0}\big[B_{\mathrm{emp}}(x,N)-B(x)\big].
\end{equation}
It is known (see e.g. \cite{Takacs 1967}, p. 170) that both \eqref{1'} and \eqref{1''}
have the same limiting distribution:
\begin{equation}\label{2'}
\begin{aligned}
F(z)&=\lim_{N\to\infty}\mathrm{Pr}\left\{\delta_N^-\leq\frac{z}{\sqrt{N}}\right\}=
\lim_{N\to\infty}\mathrm{Pr}\left\{\delta_N^+\leq\frac{z}{\sqrt{N}}\right\}\\
&=
\begin{cases}1-\mathrm{e}^{-2z^{2}}, &\mbox{for} \ z\geq0,\\
0, &\mbox{for} \ z<0.
\end{cases}
\end{aligned}
\end{equation}

For fixed $N$ the statistics $\delta_N^-$ and $\delta_N^+$ are dependent. However, as $N$ increases to infinity, they become asymptotically independent. Explanation of this fact is discussed later in Section \ref{sect4.2} (see Lemma \ref{lem3} and its proof). This fact is essentially used in our further analysis.

For the $GI/M/1/n$ queueing system, similar Kolmogorov's statistics are used for statistical analysis of stationary loss probabilities, where the probability distribution function of interarrival time is denoted by $A(x)$, empirical probability distribution of interarrival time based on $N$ observations of interarrival times is denoted by $A_{\mathrm{emp}}(x,N)$, and statistics $\delta_N$, $\delta_N^-$ and $\delta_N^+$ are, correspondingly, as follows:
$$\delta_N=\sup_{x>0}\big|A(x)-A_{\mathrm{emp}}(x,N)\big|,$$
$$
\delta_N^-=\sup_{x>0}\big[A(x)-A_{\mathrm{emp}}(x,N)\big],
$$
and
$$
\delta_N^+=\sup_{x>0}\big[A_{\mathrm{emp}}(x,N)-A(x)\big].
$$
The aim of the present paper is twofold.
Along with obtaining the ranges for aforementioned characteristics of $M/GI/1/n$ and $GI/M/1/n$
queueing systems for a given probability $P$
(which is the first aim of the paper) we also answer to the question:
\textit{Which of these Kolmogorov's statistics is better?} That is,
under which of these Kolmogorov's statistics the difference between
an upper and lower bounds based on the chosen probability $P$ is smaller?

The rest of the paper is organized as follows. In Section 2, we recall the known
representations for the characteristics studied in the paper for $M/GI/1/n$ and $GI/M/1/n$
queueing systems. In Section 3, we prove the main lemmas, which are then used to find the
estimators for required characteristics of the queueing systems. The estimators themselves
are derived in Section 4. In Section 5 we compare our statistical results and address the
question formulated above. In Section 6 numerical results are presented.
In Section 7 we conclude the paper.

\section{The recurrence relation for main characteristics of the $M/GI/1/n$ and $GI/M/1/n$ queueing systems}

In this section, we recall the known results for the main
characteristics of the $M/GI/1/n$ and $GI/M/1/n$ queueing systems
during their busy periods. For a more detailed information see
\cite{A0}, \cite{Abook} or \cite{A4}.

\subsection{The $M/GI/1/n$ queueing system} Consider the $M/GI/1/n$ queueing system, in which
the arrival flow is Poisson with parameter $\lambda$, and the
probability distribution function of the service time is $B(x)$
having the expectation $b$. Parameter $n$ denotes the number of
waiting places, i.e. the capacity for the customer in service is not
taken into account. Let $T_n$ denote the length of a busy period of
this system, and let $\nu_n$ and $L_n$ denote the number of served
and, respectively, lost customers during that busy period. The
recurrence relation for $\mathrm{E}T_n$ has been originally obtained
by Tomko \cite{Tomko}:
\begin{equation}\label{3}
\mathrm{E}T_n=\sum_{i=0}^n\mathrm{E}T_{n-i+1}\int_0^\infty\mathrm{e}^{-\lambda x}\frac{(\lambda x)^i}{i!}\mathrm{d}B(x),
\end{equation}
where $\mathrm{E}T_0=b$, and $T_i$ denotes the length of a busy period in the $M/GI/1/i$ queueing systems having the same arrival rate $\lambda$ and the same probability distribution of the service time as the original queueing system $M/GI/1/n$.
The expected number of served customers during the same busy period satisfies the recurrence
relation similar to \eqref{3}:
\begin{equation}\label{3add}
\mathrm{E}\nu_n=\sum_{i=0}^n\mathrm{E}\nu_{n-i+1}\int_0^\infty\mathrm{e}^{-\lambda x}\frac{(\lambda x)^i}{i!}\mathrm{d}B(x),
\end{equation}
where $\nu_i$ denotes the number of served customers during a busy
period $T_i$ (e.g. see \cite{A00}). The main difference between
recurrence relations \eqref{3} and \eqref{3add} is that \eqref{3add}
starts from $\mathrm{E}\nu_0=1$, while \eqref{3} starts from
$\mathrm{E}T_0=b$.

For the expected number of losses during a busy period we
correspondingly have the following recurrence relation (see e.g.
\cite{A0}, \cite{A00}, \cite{Abook} or \cite{A4} for more details),
which is similar to the previous two recurrence relations given by
\eqref{3} and \eqref{3add}:
\begin{equation}\label{4}
\mathrm{E}L_n-1=\sum_{i=0}^n(\mathrm{E}L_{n-i+1}-1)\int_0^\infty\mathrm{e}^{-\lambda x}\frac{(\lambda x)^i}{i!}\mathrm{d}B(x),
\end{equation}
where $\mathrm{E}L_0=\lambda b$, and $L_i$ denotes the number of losses during a busy period $T_i$.

All \eqref{3}, \eqref{3add} and \eqref{4} are the following convolution type recurrence relation:
\begin{equation}\label{5}
Q_n=\sum_{i=0}^n Q_{n-i+1}r_i,
\end{equation}
where $r_i=\int_0^\infty\mathrm{e}^{-\lambda x}\frac{(\lambda x)^i}{i!}\mathrm{d}B(x)$.

\subsection{The $GI/M/1/n$ queueing system} Consider the $GI/M/1/n-1$ queueing system, where an interarrival time has the probability distribution function $A(x)$, and the parameter of an exponentially distributed service time is $\mu$. (The number of waiting places $n-1$ excludes the place for a customer in service.) Let $\pi_{n}$ denote the stationary loss probability. It is shown in \cite{A01} that $\pi_n$ satisfied the recurrence relation:
\begin{equation}\label{6}
\frac{1}{\pi_n}=\sum_{i=0}^n\frac{1}{\pi_{n-i+1}}\int_0^\infty\mathrm{e}^{-\mu x}\frac{(\mu x)^i}{i!}\mathrm{d}A(x),
\end{equation}
where $\pi_0=1$. Recurrence relation \eqref{6} is similar to that of \eqref{3}, \eqref{3add} and \eqref{4}, and has general form \eqref{5}. The only difference is that in representation \eqref{6} $r_i=\int_0^\infty\mathrm{e}^{-\mu x}\frac{(\mu x)^i}{i!}\mathrm{d}A(x)$, while in representations \eqref{3}, \eqref{3add} and \eqref{4} $r_i=\int_0^\infty\mathrm{e}^{-\lambda x}\frac{(\lambda x)^i}{i!}\mathrm{d}B(x)$.

Another representation for $\pi_n$ has been obtained by Miyazawa \cite{Miyazawa}.

\section{Main lemmas} In this section we prove main lemmas, which help us to obtain then the desired estimators for characteristics of queueing systems studied in this paper.

Let $F_1(x)$ and $F_2(x)$ denote arbitrary probability distribution functions of positive random variables. For a positive parameter $\alpha$, let us denote
$$
r_i(F_1)=\int_0^\infty\mathrm{e}^{-\alpha x}\frac{(\alpha x)^i}{i!}\mathrm{d}F_1(x),
$$
and
$$
r_i(F_2)=\int_0^\infty\mathrm{e}^{-\alpha x}\frac{(\alpha x)^i}{i!}\mathrm{d}F_2(x).
$$

\begin{lem}\label{lem1}
Assume that $\sup_{x>0}\big|F_1(x)-F_2(x)\big|<\epsilon$. Then,
\begin{equation}\label{7}
\big|r_0(F_1)-r_0(F_2)\big|<\epsilon,
\end{equation}
and for all $i=1,2,...$,
\begin{equation}\label{8}
\big|r_i(F_1)-r_i(F_2)\big|<2\epsilon.
\end{equation}
\end{lem}

\begin{proof} By partial integration we have
\begin{equation}\label{8''}
\int_0^\infty\mathrm{e}^{-\alpha x}\mathrm{d}F_k(x)=\alpha\int_0^\infty\mathrm{e}^{-\alpha x}F_k(x) \mathrm{d}x, \ k=1,2.
\end{equation}
Therefore,
\begin{equation*}
\begin{aligned}
\left|\int_0^\infty\mathrm{e}^{-\alpha x}\mathrm{d}F_1(x)-\int_0^\infty\mathrm{e}^{-\alpha x}\mathrm{d}F_2(x)\right|&=\alpha\left|\int_0^\infty\mathrm{e}^{-\alpha x}F_1(x) \mathrm{d}x
-\int_0^\infty\mathrm{e}^{-\alpha x}F_2(x) \mathrm{d}x\right|\\
&\leq\underbrace{\alpha\int_0^\infty\mathrm{e}^{-\alpha x}\mathrm{d}x}_{=1}\underbrace{\left(\sup_{x>0}\big|F_1(x)-F_2(x)\big|\right)}_{<\epsilon \ \mbox{by~the~assumption}}\\
&<\epsilon.
\end{aligned}
\end{equation*}
Inequality \eqref{7} is proved.

Let us now prove inequalities \eqref{8}. For $i=1,2,...$ and $k=1,2$ by partial integration we have:
\begin{equation}\label{8'''}
\begin{aligned}
&\int_0^\infty\mathrm{e}^{-\alpha x}\frac{(\alpha x)^i}{i!}\mathrm{d}F_k(x)\\
&=\alpha\left(\int_0^\infty\mathrm{e}^{-\alpha x}\frac{(\alpha x)^i}{i!}F_k(x)\mathrm{d}x
-\int_0^\infty\mathrm{e}^{-\alpha x}\frac{(\alpha x)^{i-1}}{(i-1)!}F_k(x)\mathrm{d}x\right)
\end{aligned}
\end{equation}
Therefore,
\begin{equation}\label{9}
\begin{aligned}
&\left|\int_0^\infty \mathrm{e}^{-\alpha x}\frac{(\alpha x)^i}{i!} \mathrm{d}F_1(x)-
\int_0^\infty \mathrm{e}^{-\alpha x}\frac{(\alpha x)^i}{i!} \mathrm{d}F_2(x)\right|\\
&=\alpha\left|\int_0^\infty \mathrm{e}^{-\alpha x}\frac{(\alpha x)^i}{i!} F_1(x)\mathrm{d}x-
\int_0^\infty \mathrm{e}^{-\alpha x}\frac{(\alpha x)^i}{i!}F_2(x) \mathrm{d}x\right.\\
& \ \ \ -\left.\int_0^\infty \mathrm{e}^{-\alpha x}\frac{(\alpha x)^{i-1}}{(i-1)!} F_1(x)\mathrm{d}x+
\int_0^\infty \mathrm{e}^{-\alpha x}\frac{(\alpha x)^{i-1}}{(i-1)!}F_2(x)\mathrm{d}x\right|\\
&\leq\frac{\alpha}{i!}\int_0^\infty \mathrm{e}^{-\alpha x}(\alpha x)^i \mathrm{d}x\underbrace{\left(
\sup_{x>0}\big|F_1(x)-F_2(x)\big|\right)}_{<\epsilon \ \mbox{by~the~assumption}}\\
&\ \ \ +\frac{\alpha}{(i-1)!}\int_0^\infty \mathrm{e}^{-\alpha x}(\alpha x)^{i-1} \mathrm{d}x
\underbrace{\left(
\sup_{x>0}\big|F_1(x)-F_2(x)\big|\right)}_{<\epsilon \ \mbox{by~the~assumption}}.
\end{aligned}
\end{equation}
Notice, that
$$
\int_0^\infty \mathrm{e}^{-\alpha x}(\alpha x)^i\mathrm{d}x=\frac{1}{\alpha}\int_0^\infty
\mathrm{e}^{-y}y^i\mathrm{d}y=\frac{1}{\alpha}\Gamma(i+1),
$$
and
$$
\int_0^\infty \mathrm{e}^{-\alpha x}(\alpha
x)^{i-1}\mathrm{d}x=\frac{1}{\alpha}\int_0^\infty
\mathrm{e}^{-y}y^{i-1}\mathrm{d}y=\frac{1}{\alpha}\Gamma(i),
$$
where $\Gamma(x)$ is Euler's Gamma function. Taking into account that $\Gamma(i+1)=i!$, from the last inequality of \eqref{9} we arrived at the desired result $|r_i(F_1)-r_i(F_2)|<2\epsilon$ for $i=1,2,...$. Inequalities \eqref{8} are proved, and the proof of the lemma is completed.
\end{proof}

\begin{lem}\label{lem2}
Assume that $\sup_{x>0}\big[F_1(x)-F_2(x)\big]<\epsilon_1$ and $\sup_{x>0}\big[F_2(x)-F_1(x)\big]<\epsilon_2$. Then,
\begin{equation}\label{7'}
r_0(F_1)-r_0(F_2)<\epsilon_1,
\end{equation}
and for all $i=1,2,...$,
\begin{equation}\label{8'}
r_i(F_1)-r_i(F_2)<\epsilon_1+\epsilon_2.
\end{equation}
\end{lem}

\begin{proof} The proof of this lemma is similar to that of Lemma \ref{lem1}. From \eqref{8''} we have:
\begin{equation*}
\begin{aligned}
\int_0^\infty\mathrm{e}^{-\alpha x}\mathrm{d}F_1(x)-\int_0^\infty\mathrm{e}^{-\alpha x}\mathrm{d}F_2(x)&=\alpha\int_0^\infty\mathrm{e}^{-\alpha x}F_1(x) \mathrm{d}x
-\int_0^\infty\mathrm{e}^{-\alpha x}F_2(x) \mathrm{d}x\\
&\leq\underbrace{\alpha\int_0^\infty\mathrm{e}^{-\alpha x}\mathrm{d}x}_{=1}\underbrace{\left(\sup_{x>0}\big[F_1(x)-F_2(x)\big]\right)}_{<\epsilon_1 \ \mbox{(assumption)}}\\
&<\epsilon_1.
\end{aligned}
\end{equation*}
Inequality \eqref{7'} is proved. Next, from \eqref{8'''} we have:
\begin{equation*}\label{9'}
\begin{aligned}
&\int_0^\infty \mathrm{e}^{-\alpha x}\frac{(\alpha x)^i}{i!} \mathrm{d}F_1(x)-
\int_0^\infty \mathrm{e}^{-\alpha x}\frac{(\alpha x)^i}{i!} \mathrm{d}F_2(x)\\
&=\alpha\int_0^\infty \mathrm{e}^{-\alpha x}\frac{(\alpha x)^i}{i!} F_1(x)\mathrm{d}x-
\int_0^\infty \mathrm{e}^{-\alpha x}\frac{(\alpha x)^i}{i!}F_2(x) \mathrm{d}x\\
& \ \ \ -\int_0^\infty \mathrm{e}^{-\alpha x}\frac{(\alpha x)^{i-1}}{(i-1)!} F_1(x)\mathrm{d}x+
\int_0^\infty \mathrm{e}^{-\alpha x}\frac{(\alpha x)^{i-1}}{(i-1)!}F_2(x)\mathrm{d}x\\
&\leq\frac{\alpha}{i!}\underbrace{\int_0^\infty \mathrm{e}^{-\alpha x}(\alpha x)^i \mathrm{d}x}_{=\frac{i!}{\alpha}}\underbrace{\left(
\sup_{x>0}\big[F_1(x)-F_2(x)\big]\right)}_{<\epsilon_1 \ \mbox{(assumption)}}\\
&\ \ \ +\frac{\alpha}{(i-1)!}\underbrace{\int_0^\infty \mathrm{e}^{-\alpha x}(\alpha x)^{i-1} \mathrm{d}x}_{=\frac{(i-1)!}{\alpha}}
\underbrace{\left(
\sup_{x>0}\big[F_2(x)-F_1(x)\big]\right)}_{<\epsilon_2 \ \mbox{(assumption)}}\\
&<\epsilon_1+\epsilon_2.
\end{aligned}
\end{equation*}
Desired inequality \eqref{8'} follows. The proof is complete.
\end{proof}

\section{Explicit recursion and estimates for characteristics of queueing systems}
The recurrence relation in form \eqref{5} is not convenient for estimating the aforementioned characteristics of the queueing systems. In the left-hand side of \eqref{5} $Q_n$ is, while in the right-hand side of \eqref{5} the linear combination of $Q_{1}$, $Q_{2}$,\ldots, $Q_{n+1}$ is presented. In order to derive the appropriate estimates, we first should rewrite \eqref{5} in the form of explicit recursion. We have:
\begin{eqnarray}
Q_1&=&\frac{1}{r_0}Q_0,\label{10}\\
Q_2&=&\frac{1-r_1}{r_0}Q_1,\label{11}\\
Q_{n}&=&\frac{1}{r_0}\left[(1-r_1)Q_{n-1}-r_2Q_{n-2}-\ldots-r_{n-1}Q_1\right],
\ n\geq3.\label{12}
\end{eqnarray}
Explicit recurrence relations of \eqref{10}-\eqref{12} can be now used to establish necessary estimates for characteristics of the loss queueing systems $M/GI/1/n$ and $GI/M/1/n$.
To build the recursion for $\mathrm{E}T_n$, we start from the $M/GI/1/0$ queueing system. The busy period of this system contains only a single service time. In the sequel we assume that the parameters $\lambda$ and $b$ of the queueing systems are known. As the expected service time $b$ is known, we set $\widehat T_0=b$, that is, the expected busy period of the $M/GI/1/0$ queueing system is reckoned to be estimated exactly. In the following we use the notation:
$$
r_i(B)=\int_0^\infty\mathrm{e}^{-\lambda x}\frac{(\lambda x)^i}{i!}\mathrm{d}B(x),
$$
and
$$
r_i(B_{\mathrm{emp}})=\int_0^\infty\mathrm{e}^{-\lambda x}\frac{(\lambda x)^i}{i!}\mathrm{d}B_{\mathrm{emp}}(x,N),
$$
where $B_{\mathrm{emp}}(x,N)$ is an empirical probability distribution function based on $N$ observations. In the following the estimator for $Q_n$ will be denoted $\widehat Q_n$ and specifically that for the expected busy period will be denoted $\widehat T_n$.
(Similarly, for $GI/M/1/n$ queueing systems, the notation $r_i(A)$ and $r_i(A_{\mathrm{emp}})$ can be used where the parameter $\lambda$ should be replaced with the parameter $\mu$.)

\subsection{Estimators based on the statistic $\delta_N$}\label{sect4.1}
Assume that the inequality
\begin{equation}\label{13}
\delta_N<\epsilon
\end{equation}
holds with probability $P$. Since $N$ is assumed to be large enough, this probability can be chosen from limit relation \eqref{2}.

Then, according to Lemma \ref{lem1}, with the probability not smaller than $P$ we have
$$
\big|r_0(B)-r_0(B_{\mathrm{emp}})\big|<\epsilon,
$$
and for $i=1,2,\ldots,n$
$$
\big|r_i(B)-r_i(B_{\mathrm{emp}})\big|<2\epsilon.
$$
For further simplifications, we will write $r_i=r_i(B_{\mathrm{emp}})$, $i=0,1,2,\ldots$, omitting the argument $B_{\mathrm{emp}}$ in the notation.

On the basis \eqref{10}-\eqref{12} and Lemma \ref{lem1} we have the following.
\begin{thm}\label{thm-1}
The point estimator $\widehat Q_n$ for a required characteristic of
a queueing system is recurrently defined as
\begin{eqnarray}
\widehat{Q}_1&=&\frac{1}{r_0}Q_0,\label{10'}\\
\widehat{Q}_2&=&\frac{1-r_1}{r_0}\widehat{Q}_1,\label{11'}\\
\widehat{Q}_{n}&=&\frac{1}{r_0}\left[(1-r_1)\widehat
Q_{n-1}-r_2\widehat Q_{n-2}-\ldots-r_{n-1}\widehat Q_1\right], \
n\geq3.\label{12'}
\end{eqnarray}
Then, the interval estimators with a confidence probability
non-smaller than $P$ are recurrently defined as
\begin{equation}\label{14}
\widehat Q_1^{\mathrm{lower}}=\frac{1}{r_0+\epsilon}Q_0<
Q_1<\frac{1}{r_0-\epsilon}Q_0=\widehat Q_1^{\mathrm{upper}},
\end{equation}
\begin{equation}\label{15}
\widehat
Q_2^{\mathrm{lower}}=\frac{1-r_1-2\epsilon}{[r_0+\epsilon]^{2}}Q_0 <
 Q_2<\frac{1-r_1+2\epsilon}{[r_0-\epsilon]^{2}}Q_0 =
Q_2^{\mathrm{upper}},
\end{equation}
\begin{equation}\label{16}
\begin{aligned}
\widehat Q_{n}^{\mathrm{lower}}&=\frac{1}{r_0+\epsilon}\left[
(1-r_1-2\epsilon)\widehat Q_{n-1}^{\mathrm{lower}}\right.\\
&\ \ \ \left.-(r_2+2\epsilon) \widehat Q_{n-2}^{\mathrm{upper}}
-\ldots-(r_{n-1}+2\epsilon)\widehat Q_{1}^{\mathrm{upper}}
\right]\\
&< Q_{n}<\frac{1}{r_0-\epsilon}\left[
(1-r_1+2\epsilon)\widehat Q_{n-1}^{\mathrm{upper}}\right.\\
&\ \ \ \left.-(r_2-2\epsilon) \widehat Q_{n-2}^{\mathrm{lower}}
-\ldots-(r_{n-1}-2\epsilon)\widehat Q_{1}^{\mathrm{lower}}
\right]\\
&=\widehat Q_{n}^{\mathrm{upper}}, \ n\geq3.
\end{aligned}
\end{equation}
The value $\epsilon$ is determined from the equation:
$K\left({\epsilon}{\sqrt{N}}\right)=P$, $K(z)$ is Kolmogorov's
function defined in \eqref{2}, and in the notation of
\eqref{14}-\eqref{16} $\widehat Q_{i}^{\mathrm{lower}}$ and
$\widehat Q_{i}^{\mathrm{upper}}$ are used for lower and,
respectively, upper bound of the interval estimators.
\end{thm}

\begin{proof} Indeed for any $0<\epsilon<r_0$ we apparently
have
\begin{equation}\label{eq.add.1}
\widehat
{Q}_1^{\mathrm{lower}}=\frac{1}{r_0+\epsilon}Q_0<\frac{1}{r_0}Q_0<\frac{1}{r_0-\epsilon}Q_0=\widehat
{Q}_1^{\mathrm{upper}},
\end{equation}
where $\frac{1}{r_0}Q_0$ in the middle of inequality
\eqref{eq.add.1} is $\widehat Q_1$. This means that for $\epsilon>0$
given such that $|r_0(B)-r_0(B_{\mathrm{emp}}|<\epsilon$ occurs with
confidence probability not smaller than $P$, we have
$$
\widehat
{Q}_1^{\mathrm{lower}}=\frac{1}{r_0+\epsilon}Q_0<Q_1<\frac{1}{r_0-\epsilon}Q_0=\widehat
{Q}_1^{\mathrm{upper}}.
$$
Relation \eqref{15} follows similarly. Indeed, for any
$0<\epsilon<r_0$ we apparently have
\begin{equation}\label{eq.add.2}
\frac{1-r_1-2\epsilon}{[r_0+\epsilon]^{2}}Q_0
<\frac{1-r_1}{r_0^{2}}Q_0
<\frac{1-r_1+2\epsilon}{[r_0-\epsilon]^{2}}Q_0,
\end{equation}
where $\frac{1-r_1}{r_0^{2}}Q_0$ in the middle of inequality
\eqref{eq.add.2} is $\widehat{Q}_2$. That is, for any $\epsilon>0$
given such that $|r_0(B)-r_0(B_{\mathrm{emp}}|<\epsilon$ occurs with
the confidence probability not smaller than $P$, we have
$$
\frac{1-r_1-2\epsilon}{[r_0+\epsilon]^{2}}Q_0 <Q_2
<\frac{1-r_1+2\epsilon}{[r_0-\epsilon]^{2}}Q_0.
$$
Then \eqref{16} is proved by the similar arguments with the aid of
induction. Namely, for $n\geq3$ and any $0<\epsilon<r_0$ we have the
inequality
\begin{equation*}
\begin{aligned}
&\frac{1}{r_0+\epsilon}\left[ (1-r_1-2\epsilon)\widehat
Q_{n-1}^{\mathrm{lower}}-(r_2+2\epsilon) \widehat
Q_{n-2}^{\mathrm{upper}} -\ldots-(r_{n-1}+2\epsilon)\widehat
Q_{1}^{\mathrm{upper}} \right]< \widehat
Q_{n}\\&<\frac{1}{r_0-\epsilon}\left[ (1-r_1+2\epsilon)\widehat
Q_{n-1}^{\mathrm{upper}}-(r_2-2\epsilon) \widehat
Q_{n-2}^{\mathrm{lower}}-\ldots-(r_{n-1}-2\epsilon)\widehat
Q_{1}^{\mathrm{lower}} \right],
\end{aligned}
\end{equation*}
which means that for $\epsilon>0$ given such that
$|r_0(B)-r_0(B_{\mathrm{emp}}|<\epsilon$ occurs with the confidence
probability not smaller than $P$, we have
\begin{equation*}
\begin{aligned}
&\frac{1}{r_0+\epsilon}\left[ (1-r_1-2\epsilon)\widehat
Q_{n-1}^{\mathrm{lower}}-(r_2+2\epsilon) \widehat
Q_{n-2}^{\mathrm{upper}} -\ldots-(r_{n-1}+2\epsilon)\widehat
Q_{1}^{\mathrm{upper}} \right]<
Q_{n}\\&<\frac{1}{r_0-\epsilon}\left[ (1-r_1+2\epsilon)\widehat
Q_{n-1}^{\mathrm{upper}}-(r_2-2\epsilon) \widehat
Q_{n-2}^{\mathrm{lower}}-\ldots-(r_{n-1}-2\epsilon)\widehat
Q_{1}^{\mathrm{lower}} \right].
\end{aligned}
\end{equation*}
Theorem \ref{thm-1} is proved.
\end{proof}

Using relations \eqref{10'}-\eqref{16} one can build estimators and
the bounds for these estimators for all of the aforementioned
characteristics of the queueing systems $M/GI/1/n$ and $GI/M/1/n$.

\begin{rem}
For relations \eqref{14}-\eqref{16} the following conventions are
made. If $r_0\leq\epsilon$, then the right-hand sides of
\eqref{14}-\eqref{16} are set to be equal to infinity. If
$1-r_1-2\epsilon$ is negative, then the left-hand side of \eqref{15}
is set to be equal to zero. In addition, if one or other term
$1-r_i-2\epsilon$, $i=1,2,\ldots,n$ of the left-hand side of
\eqref{16} is negative, that term is set to zero. The value
$\epsilon$ should be chosen small enough, i.e. appropriate sample
size $N$ should be chosen large enough, such that all of the above
inconsistences can occur with negligibly small probability.
\end{rem}

\subsection{Estimators based on the statistic $\delta_N^-$ or $\delta_N^+$}\label{sect4.2}
Since both of these statistics are similar, we only consider the
first one, $\delta_N^-$. Assume that $N$ is large enough and
$\delta_N^-<\epsilon$ occurs with the probability not smaller than
$P$.

According to Lemma 3.2 with the same probability $P$
$$
r_0(B)-r_0(B_{\mathrm{emp}})<\epsilon,
$$
and
$$
r_0(B_{\mathrm{emp}})-r_0(B)<\epsilon.
$$
From same Lemma \ref{lem2} according to \eqref{8'}  we have
$r_i(B)-r_i(B_{\mathrm{emp}})<\delta_N^-+\delta_N^+$ for
$i=1,2,\ldots,n$. Therefore, given that
$\mathrm{Pr}\{\delta_N^{+}+\delta_N^{-}<\gamma\}=P$, with the
confidence probability not smaller than $P$, we obtain the
inequalities:
$$
r_i(B)-r_i(B_{\mathrm{emp}})<\gamma,
$$
and
$$
r_i(B_{\mathrm{emp}})-r_i(B)<\gamma.
$$

\begin{lem}\label{lem3}
As $N$ increases to infinity, the statistics $\sqrt{N}\delta_N^{+}$ and $\sqrt{N}\delta_N^{-}$ become asymptotically independent.
\end{lem}

\begin{proof}
This fact can be proved with the aid of construction used in \cite{Takacs 1967}, p.171-172 as follows.
Assume that $B(x)$ is a continuous probability distribution function, and $B_{\mathrm{emp}}(x,N)$ is the empirical probability distribution function based on the sample ($\xi_1$, $\xi_2$,\ldots,$\xi_N$), that is, $B_{\mathrm{emp}}(x,N)$ is defined as the fraction of random variables that less than or equal to $x$. Let us write the sample ($\xi_1$, $\xi_2$,\ldots,$\xi_N$) in ascending order as
($\xi_1^*$, $\xi_2^*$,\ldots,$\xi_N^*$). Consider the deviations $\delta_N^-(n)=B(\xi_n^*)-B_{\mathrm{emp}}(\xi_n^*,N)$ and
$\delta_N^+(n)=B_{\mathrm{emp}}(\xi_n^*,N)-B(\xi_n^*)$, $n=1,2,\ldots,N$. Apparently, the random variables $\delta_N^-(n)$ and $\delta_N^+(n)$ all are continuous ($n=1,2,\ldots,N$), and all their joint distributions are independent of the choice of the probability distribution $B(x)$. That is, without loss of generality one can assume that $B(x)$ is the uniform distribution in [0,1]. Then, $\delta_N^+(n)= \xi_n^*-\frac{n}{N}$ and $\delta_N^-(n)= \frac{n}{N}-\xi_n^*$. Let $\kappa_n$ be the number of random variables amongst $\xi_1$, $\xi_2$,\ldots, $\xi_N$ falling into the interval $\left(\frac{n-1}{N},\frac{n}{N}\right]$ and $\mathcal{N}_n=\kappa_1+\kappa_2+\ldots+\kappa_n$, $n=1,2,\ldots,N$. Then the number of points where $\delta_N^-(n)$ are positive is $\sum_{n=1}^N\mathrm{I}\{\mathcal{N}_n\geq n\}$, and, respectively, the number of points where $\delta_N^+(n)$ are positive is $\sum_{n=1}^N\mathrm{I}\{\mathcal{N}_n< n\}$. For $N$ fixed, the random variables $\delta_N^-=\max_{1\leq n\leq N} \delta_N^-(n)$ and $\delta_N^+=\max_{1\leq n\leq N} \delta_N^+(n)$ are generally dependent. As $N\to\infty$, the random variables $\delta_N^+$ and $\delta_N^-$ both vanish. The aim is to prove that the normalized random variables $\sqrt{N} \delta_N^+$ and $\sqrt{N} \delta_N^-$ become asymptotically independent as $N\to\infty$.

For this purpose we will use the known fact that the random variable $\xi_n^*$ has the Beta distribution
with parameters $n$ and $N-n+1$ (see e.g. David and Nagaraja \cite{DN}, p.14).
In other words, this means that the random variable $\xi_n^*$ has the same distribution as
$$
\frac{\eta_1+\eta_2+\ldots+\eta_n}{\eta_1+\eta_2+\ldots+\eta_{N+1}},
$$
where $\eta_1$, $\eta_2$,\ldots, $\eta_N$,\ldots is the sequence of independent exponentially
distributed random variables with parameter 1 (see Karlin \cite{Kar}, p.242-244).
According to the strong law of large numbers,
$$
\mathrm{Pr}\left\{\lim_{N\to\infty}\frac{N}{\eta_1+\eta_2+\ldots+\eta_{N+1}}=1\right\}=1.
$$
Therefore, as $N\to\infty$, the random variable $N\delta_N^+(n)=N\xi_n^*-n$ can be asymptotically
represented as
$$
(\eta_1-1)+(\eta_2-1)+\ldots+(\eta_n-1).
$$
From this asymptotic representation, apparently that $\lim_{N\to\infty}\mathrm{E}(N\xi_n^*-n)=0$ for all
$n=1,2,\ldots,N$, and the expected number of positive and negative values of $N\xi_n^*-n$
for $n=1,2,\ldots,N$ is approximately the same. That is the numbers of points where $\delta_N^+(n)$ and
$\delta_N^-(n)$ are positive are asymptotically equal.

For large $N$, let $j$ denote an index such that $N\delta_N^+(j)$ is negative while $N\delta_N^+(j+1)$
is positive. Then there is a random number $\upsilon$ of consequent positive values
of $N\delta_N^+(j+i)$, $i=1,2,\ldots,\upsilon$. This random number is called
\textit{length of positive period} or simply \textit{positive period}. Notice that,
\begin{equation}\label{rrel}
\mathrm{Pr}\left\{\sum_{i=1}^{j+1}(\eta_i-1)\leq x\Big|\sum_{i=1}^{j}
(\eta_i-1)<0, \sum_{i=1}^{j+1}(\eta_i-1)>0\right\}=1-\mathrm{e}^{-x}.
\end{equation}
Indeed, let $\tau_j=j-\sum_{i=1}^j\eta_i+1$. Then, the conditional probability of the
left-hand side of \eqref{rrel} can be written
\begin{equation*}
\begin{aligned}
\mathrm{Pr}\left\{\eta_{j+1}\leq x+\tau_j|\eta_{j+1}>\tau_j>1\right\}&=\int_0^\infty \mathrm{Pr}\left\{\eta_{j+1}\leq x+y|\eta_{j+1}>y\right\}\mathrm{d}\mathrm{Pr}\{\tau_j\leq y\}\\
&=(1-\mathrm{e}^{-x})\int_0^\infty\mathrm{d}\mathrm{Pr}\{\tau_j\leq y\}\\
&=1-\mathrm{e}^{-x}.
\end{aligned}
\end{equation*}
\eqref{rrel} follows. Thus, the limiting distribution of $N\delta_N^+(j+1)$ as $N\to\infty$
is exponential, and $N\delta_N^+(j+1)$ is asymptotically independent of the past
values of $N\delta_N^+(k)$, $k=1,2,\ldots, j$. The following random variables
$N\delta_N^+(j+i)$, $i=2,3,\ldots,\upsilon$ are also asymptotically independent of the
aforementioned past values $N\delta_N^+(k)$, $k=1,2,\ldots, j$. So, the positive period and
associated random variables $N\delta_N^+(j+i)$, $i=2,3,\ldots,\upsilon$ all are asymptotically
independent of the past.

Consequently, let $j_1$, $j_2$,\ldots be the sequence of values such that $N\delta_N^+(j_k)$,
$k=1,2,\ldots$ are negative while $N\delta_N^+(j_k+1)$ are positive.
Then, according to the above statement, the positive periods are asymptotically mutually
independent and identically distributed random variables.

The similar result can be obtained for \textit{negative} periods
which are constructed similarly to those positive periods. Indeed,
for $N$ large, let now $j$ denote an index such that
$N\delta_N^+(j)$ is positive while $N\delta_N^+(j+1)$ is negative.
Then there is a random number $\iota$ of consequent negative values
of $N\delta_N^+(j+i)$, $i=1,2,\ldots,\iota$. This random number is
called \textit{length of negative period} or simply \textit{negative
period}. Notice that the event $\{N\delta_N^+(j)>0 \ \text{and} \
N\delta_N^+(j+1)<0\}$ means that during the time interval $[j,j+1)$
there are at least two events of Poisson process with rate 1, and
the second event can be considered as the start of a negative
period. As in the case of positive period, the distribution of the
length of a negative period does not depend on the past history as
well.

Thus, as $N\to\infty$, positive and negative periods, alternatively
changing one another, are asymptotically mutually independent and
identically distributed random variables.

As $N\to\infty$, the sets of values  where $N\delta_N^+(n)$ and $N\delta_N^-(n)$ are positive, are asymptotically independent, because the positive and negative periods are asymptotically independent. Note, that the set of indexes where the periods are positive is $\{n\leq N: \mathcal{N}_n<n\}$ and that set where periods are negative is $\{n\leq N: \mathcal{N}_n\geq n\}$.
Furthermore,
$$
\mathrm{Pr}\left\{\lim_{N\to\infty}\frac{1}{N}\sum_{n=1}^N\mathrm{I}\{\mathcal{N}_n<n\}
=\frac{1}{2}\right\}=1
$$

On the other hand, as $N\to\infty$, the limiting distributions of $\sqrt{N}\delta_N^-$ and $\sqrt{N}\delta_N^+$ do exist, and they are defined by \eqref{2'}.

Therefore,
\begin{eqnarray*}
\begin{aligned}
&\lim_{N\to\infty}\mathrm{Pr}
\left\{\sqrt{N}\delta_N^+\leq x, \sqrt{N}\delta_N^-\leq y \right\}\\
&=\lim_{N\to\infty}\mathrm{Pr}\left\{\frac{1}{\sqrt{N}}\max_{\{n\leq N: \mathcal{N}_n<n\}}N\delta_N^+(n)\leq x, \frac{1}{\sqrt{N}}\max_{\{n\leq N: \mathcal{N}_n\geq n\}}N\delta_N^-(n)\leq y\right\}\\
&=\lim_{N\to\infty}\mathrm{Pr}\left\{\frac{1}{\sqrt{N}}\max_{\{n\leq N: \mathcal{N}_n<n\}}N\delta_N^+(n)\leq x\right\}\\
&\ \ \ \times\lim_{N\to\infty}\mathrm{Pr}\left\{\frac{1}{\sqrt{N}}\max_{\{n\leq N: \mathcal{N}_n\geq n\}}N\delta_N^-(n)\leq y\right\}\\
&=\lim_{N\to\infty}\mathrm{Pr}
\left\{\sqrt{N}\delta_N^+\leq x\right\}\lim_{N\to\infty}\mathrm{Pr}
\left\{\sqrt{N}\delta_N^-\leq y\right\},
\end{aligned}
\end{eqnarray*}
and the proof of the statement of the lemma is completed.
\end{proof}

From Lemma \ref{lem3} and \eqref{2'} we have:
\begin{equation}\label{3'}
\lim_{N\to\infty}\mathrm{Pr}\left\{\delta_N^-+\delta_N^+\leq\frac{z}{\sqrt{N}}\right\}=F*F(z)
\end{equation}
where the probability distribution function $F(z)$ is defined in
\eqref{2'}, and the asterisk denotes convolution of the probability
distribution function $F(z)$ with itself.

For the convolution of the probability distribution $F(z)=1-\mathrm{e}^{-2z^{2}}$ with itself we have:
\begin{equation}\label{3''}
\begin{aligned}
&\int_0^z(1-\mathrm{e}^{-2(z-x)^{2}})4x\mathrm{e}^{-2x^{2}}\mathrm{d}x
=1-\mathrm{e}^{-2z^{2}}-
\sqrt{\pi}z\mathrm{e}^{-z^{2}}\big[2\Phi(\sqrt{2}z)-1\big],
\end{aligned}
\end{equation}
where $\Phi(z)=\frac{1}{\sqrt{2\pi}}\int_{-\infty}^{z}\mathrm{e}^{-\frac{y^{2}}{2}}\mathrm{d}y.$ So, the right-hand side of \eqref{3'} is explicitly determined.

Now, we have the following theorem.
\begin{thm}\label{thm2}
The point estimator $\widehat Q_n$ for a required characteristic of
a queueing system is recurrently defined by \eqref{10'}-\eqref{12'}.
Then, the interval estimators with a confidence probability
non-smaller than $P$ are recurrently defined as
\begin{equation}\label{14'}
\widehat
Q_1^{\mathrm{lower}}=\frac{1}{r_0+\epsilon}Q_0<Q_1<\frac{1}{r_0-\epsilon}Q_0=\widehat
Q_1^{\mathrm{upper}},
\end{equation}
\begin{equation}\label{15'}
\widehat
Q_2^{\mathrm{lower}}=\frac{1-r_1-\gamma}{[r_0+\epsilon]^{2}}Q_0
<Q_2<\frac{1-r_1+\gamma}{[r_0-\epsilon]^{2}}Q_0 =\widehat
Q_2^{\mathrm{upper}},
\end{equation}
\begin{equation}\label{16'}
\begin{aligned}
\widehat Q_{n}^{\mathrm{lower}}&=\frac{1}{r_0+\epsilon}\left[
(1-r_1-\gamma)\widehat Q_{n-1}^{\mathrm{lower}}\right.\\
&\ \ \ \left.-(r_2+\gamma) \widehat Q_{n-2}^{\mathrm{upper}}
-\ldots-(r_{n-1}+\gamma)\widehat Q_{1}^{\mathrm{upper}}
\right]\\
&<Q_{n}<\frac{1}{r_0-\epsilon}\left[
(1-r_1+\gamma)\widehat Q_{n-1}^{\mathrm{upper}}\right.\\
&\ \ \ \left.-(r_2-\gamma) \widehat Q_{n-2}^{\mathrm{lower}}
-\ldots-(r_{n-1}-\gamma)\widehat Q_{1}^{\mathrm{lower}}
\right]\\
&=\widehat Q_{n}^{\mathrm{upper}}, \ n\geq3,
\end{aligned}
\end{equation}
where the value $\epsilon$ is determined from the equation: $1-\exp\left(-{2{{N}}\epsilon^2}\right)=P$,
and the value $\gamma$ is determined from the equation:
$$
1-\exp\left(-2{N}\gamma^{2}\right)- \sqrt{\pi
N}\gamma\exp\left(-{{N}\gamma^{2}}\right)
\left[2\Phi\left(\sqrt{2N}{\gamma}\right)-1\right]=P,
$$
where $\Phi(x)=\frac{1}{\sqrt{2\pi}}\int_{-\infty}^x\mathrm{e}^{-\frac{y^{2}}{2}}\mathrm{d}x$.
\end{thm}

\begin{proof}
The proof of \eqref{14'}-\eqref{16'} is similar to that proof of
\eqref{14}-\eqref{16} in Theorem \ref{thm-1}.
\end{proof}

\section{Which of the statistics is better?}
In this section we address the question which of the statistics $\delta_N$, $\delta_N^+$ or $\delta_N^-$ is better, or under which of them the difference between upper and lower bound is smaller? The statistics $\delta_N^+$ and $\delta_N^-$ are symmetric and they give the same bounds. Therefore the question of comparison should be addressed to the statistic $\delta_N$ and one of the statistics $\delta_N^+$ or $\delta_N^-$, say $\delta_N^-$.

Clearly, that
\begin{equation}\label{20}
\begin{aligned}
\delta_N&:=\sup_{x\geq0}\big|B(x)-B_{\mathrm{emp}}(x,N)\big|\\
&\geq \sup_{x\geq0}\big[B(x)-B_{\mathrm{emp}}(x,N)\big]\\
&:=\delta_N^-.
\end{aligned}
\end{equation}
Therefore, for the same probability $P$ the $z$-value of the equation $P=K(z)$ is not smaller that the $z$-value of the equation $P=1-\mathrm{e}^{-2z^2}$ for any given probability $P$.
As a result, $\epsilon(\delta_N)\geq\epsilon(\delta_N^-)$, where $\epsilon(\delta_N)$ is the value $\epsilon$ obtained for the statistic $\delta_N$, and $\epsilon(\delta_N^-)$ is the value $\epsilon$ obtained for the statistic $\delta_N^-$ both for the same value of probability $P$.

On the other hand, for any probability distribution function $G(x)$
of a positive random variable we have as follows. Let
$F_1(x)=G\left(\frac{x}{2}\right)$ and $F_2(x)=(G*G)(x)$, where the
asterisk denotes convolution. Apparently,
\begin{equation}\label{conv1}
\int_0^\infty x\mathrm{d}F_1(x)=\int_0^\infty x\mathrm{d}F_2(x),
\end{equation}
and
\begin{equation}\label{conv2}
\int_0^\infty x^2\mathrm{d}F_1(x)\geq\int_0^\infty
x^2\mathrm{d}F_2(x)
\end{equation}
(the integrals in \eqref{conv1} and \eqref{conv2} are assumed to
converge). Then (see e.g. \cite{Stoyan}) $F_2(x)$ is said to be
smaller than $F_1(x)$ in the convex sense.

Then, there exists a point $x_0$ where the probability distribution
functions $F_1(x)$ and $F_2(x)$ cut one another (see Karlin and
Novikoff \cite{KarNov} or Stoyan \cite{Stoyan}, p. 13).

In the given case $F(x)=1-\mathrm{e}^{-2x^2}$. Therefore,
$F_1(x)=1-\mathrm{e}^{-\frac{1}{2}x^2}$, while
$F_2(x)=1-\mathrm{e}^{-2x^{2}}-
\sqrt{\pi}x\mathrm{e}^{-x^{2}}\big[2\Phi(\sqrt{2}x)-1\big]$.
Apparently, there exists a point $x_0$ such that $F_1(z)\leq F_2(z)$
for all $z\geq x_0$. Numerical calculations show that
$x_0\approx1.385$, which corresponds to the level of probability
$F_1(1.385)\approx F_2(1.385)\approx0.6166.$

Let $z_1$ and $z_2$ be two $z$-points associated with the solutions
of the corresponding equations $F_1(z)=P$ and $F_2(z)=P$, where $P$
is some given level of probability. Clearly, that if both $z_1\geq
x_0$ and $z_2\geq x_0$, then, since $F_1(z)\leq F_2(z)$ in this set,
we also have $z_1\geq z_2$.

This enables us to conclude, that correspondingly to  these
$z$-points of the aforementioned probability distribution functions
we obtain $2\epsilon(\delta_N^-)\geq \gamma(\delta_N^-)$. The
following notation is used here. As before, by
$\epsilon(\delta_N^-)$ we mean the value $\epsilon$ associated with
the solution $1-\mathrm{e}^{-2N\epsilon^2}=P$. By
$\gamma(\delta_N^-)$ we mean the value $\gamma$ associated with the
solution
$$
1-\exp\left(-2{N}\gamma^{2}\right)- \sqrt{\pi
N}\gamma\exp\left(-{{N}\gamma^{2}}\right)
\left[2\Phi\left(\sqrt{2N}{\gamma}\right)-1\right]=P.
$$
The aforementioned relation $2\epsilon(\delta_N^-)\geq
\gamma(\delta_N^-)$ follows from \eqref{conv1} and \eqref{conv2} as
follows. For a given $z$ value associated with large level of
probability $P$ we have $F_1(z)\leq F_2(z)$. Hence, as it was
mentioned before, from $F_1(z_1)=P$ and $F_2(z_2)=P$ we obtain
$z_2\geq z_1$. That is, from \eqref{20} we have
$\epsilon(\delta_N^-)\leq\epsilon(\delta_N)$ and, for large enough
$z$-values taking into account that $F(z)=F_1(2z)$, we have
$\gamma(\delta_N^-)\leq2\epsilon(\delta_N^-)$, and since
$\epsilon(\delta_N^-)\leq\epsilon(\delta_N)$ we finally have
$\gamma(\delta_N^-)\leq 2\epsilon(\delta_N)$. (Recall that by
$\epsilon(\delta_N)$ we mean the value of $\epsilon$ associated with
the solution $K(\epsilon\sqrt{N})=P$.)

Thus, the statistic $\delta_N^-$ becomes better than the statistic $\delta_N$ for large enough $z$-values, i.e. for large enough values of probability $P$. This conclusion is supported by the numerical examples in the next section.

\section{Numerical examples} Numerical examples are provided for $N=10,000$. For simplicity we take $B(x)=1-\mathrm{e}^{-\mu x}$, and $\lambda=\mu=1$. The value $n$ is taken 4. In the numerical examples below we build estimators for the expected busy periods of the $M/GI/1/n$ queueing systems.

\subsection{The statistic $\delta_N$}

From the equation $K(z)=P=0.95$ we obtain $z\approx1.3581$. Therefore,
$$\epsilon\approx\frac{1.3581}{\sqrt{10,000}}=0.013581.$$

With the aid of simulation, we build empirical probability distribution,
and on the basis of that empirical probability distribution we obtained
$r_i(B_{\mathrm{emp}})$, $i=0,1,2,3,4$, given in Table 1.

\begin{table}
    \begin{center}
        \begin{tabular}{c||l|l}\hline
            &  Expected     & The value obtained\\
 Statistics & (theoretical) & by simulating     \\
            &  value        & empirical distr.  \\
            & $r_i(B)$      &$r_i(B_{\mathrm{emp}})$\\
\hline $r_0$ & .5   & .5031 \\
$r_1$        & .25   & .2488  \\
$r_2$        & .125   & .1234  \\
$r_3$        & .0625   & .0615  \\
$r_4$        & .03125   & .0308  \\
\hline
        \end{tabular}

        \caption{The table for statistics $r_i(B_{\mathrm{emp}})$ obtained by simulating
        empirical probability distribution function}
    \end{center}
\end{table}

Then, on the basis of statistics $r_i(B_{\mathrm{emp}})$, $i=0,1,2,3,4$ the lower and upper bounds
for estimators $\widehat T_i$, $i=1,2,3,4$ and these estimators themselves are shown in Table 2.

\begin{table}
    \begin{center}
        \begin{tabular}{c||l|l|l|l}\hline
            &  Expected     & Estimated & Lower & Upper\\
 Statistics & (theoretical) & value     & bound & bound\\
            &  value        &           &       &      \\
            &~~$\mathrm{E}T_i$&$~~\widehat{T}_i$&~~$\widehat{T}_i^{\mathrm{lower}}$&~~$\widehat{T}_i^{\mathrm{upper}}$\\
\hline
$T_0$        & 1   & 1  &   1     & 1\\
$T_1$        & 2   & 1.987589  & 1.935434 & 2.042817\\
$T_2$        & 3   & 2.967558  & 2.71285  & 3.248177\\
$T_3$        & 4   & 3.943322  & 3.206328 & 4.783615\\
$T_4$        & 5   & 4.916821  & 3.317455 & 7.057548\\
\hline
        \end{tabular}

        \caption{Estimators for expected busy periods in $M/GI/1/n$ queueing systems based on the
        statistics $\delta_N$}
    \end{center}
\end{table}

\subsection{The statistic $\delta_N^-$ or $\delta_N^+$} From the equation $1-\mathrm{e}^{-2z^{2}}=P=0.95$ we
obtain $z=1.224$. Therefore $\epsilon=0.01224$. From the other equation
$$1-\mathrm{e}^{-2z^{2}}-
\sqrt{\pi}z\mathrm{e}^{-z^{2}}\big[2\Phi(\sqrt{2}z)-1\big]=P=0.95$$
(see \eqref{3''}), we have $z=2.08$. Therefore $\gamma=0.0208$. In
this case the lower and upper bounds for estimators $\widehat T_i$,
$i=1,2,3,4$ are shown in Table 3. (For convenience the estimators
$\widehat T_i$, $i=1,2,3,4$ themselves are duplicated from Table 2.)
\begin{table}
    \begin{center}
        \begin{tabular}{c||l|l|l|l}\hline
            &  Expected     & Estimated & Lower & Upper\\
 Statistics & (theoretical) & value     & bound & bound\\
            &  value        &           &       &      \\
            &~~$\mathrm{E}T_i$&$~~\widehat{T}_i$&~~$\widehat{T}_i^{\mathrm{lower}}$&~~$\widehat{T}_i^{\mathrm{upper}}$\\
\hline
$T_0$        & 1   & 1  &   1     & 1\\
$T_1$        & 2   & 1.987589  & 1.940466 & 2.037241\\
$T_2$        & 3   & 2.967558  & 2.750256 & 3.204070\\
$T_3$        & 4   & 3.943322  & 3.327933 & 4.633603\\
$T_4$        & 5   & 4.916821  & 3.616202 & 6.673106\\
\hline
        \end{tabular}

        \caption{Estimators for expected busy periods in $M/GI/1/n$ queueing systems based on the statistics $\delta_N^-$ or $\delta_N^+$}
    \end{center}
\end{table}

\section{Discussion and the future work}
We provided statistical analysis of certain output characteristics of $M/GI/1/n$ and $GI/M/1/n$
queueing systems by using the known statistics associated with empirical distribution.
It follows from our results and numerical analysis that in certain cases the use of the statistic
$\delta_N^-$ or $\delta_N^+$ is more profitable than that of $\delta_N$. Namely, in these certain
cases the differences between the upper and lower bounds for the output characteristics are smaller
when we use the statistic $\delta_N^-$ or $\delta_N^+$ compared to these differences between the bounds
when we use the statistic $\delta_N$.

It should be noted, that presented methods can be used for other
queueing models, in which characteristics can be expressed via
convolution type recurrence relation \eqref{5}. An example is the
extended $M/GI/1/n$ loss queueing system describing models of
telecommunication systems considered in \cite{Abramov 2004}. Another
example is the state-depending queueing system describing models of
dam/inventory systems considered in \cite{Abramov 2007}. More
complicated models of loss queueing system with batch
arrivals/departures of customers are of special interest as well and
can be considered as a subject of future work. It is also
interesting to study the cases where the parameters $\lambda$ and/or
$b$ are unknown. This case is more realistic for practical
applications and leads to new challenging aspects of this theory.

\section*{Acknowledgements} The results of this paper are discussed in a local
seminar in Monash University (Australia). The author thanks the
colleagues and especially Prof. Kais Hamza for interesting questions
leading to a substantial improvement of the original version of the
paper. Many thanks are to the anonymous reviewer of this paper,
whose recommendations improved this paper substantially. The
financial support of the Australian Research Council is highly
appreciated as well.

\section*{Appendix: Calculation of the integral in \eqref{3''}}

We have:
$$
\begin{aligned}
\int_0^z(1-\mathrm{e}^{-2(z-x)^2})4x\mathrm{e}^{-2x^2}\mathrm{d}x&=\int_0^z4x\mathrm{e}^{-2x^2}\mathrm{d}x
-\int_0^z4x\mathrm{e}^{-2x^2}\mathrm{e}^{-2(z-x)^2}\mathrm{d}x\\
&=I_1+I_2.
\end{aligned}
$$
Clearly, that $I_1=F(z)=1-\mathrm{e}^{-2z^2}$. For $I_2$ we have the
following:
$$
\begin{aligned}
I_2&=-\int_0^z4x\mathrm{e}^{-2z^2+4zx-4x^2}\mathrm{d}x=
-\mathrm{e}^{-z^2}\int_0^z4x\mathrm{e}^{-(z-2x)^2}\mathrm{d}x\\
&=-\mathrm{e}^{-z^2}\int_0^z[(2x-z)+z]\mathrm{e}^{-(z-2x)^2}\mathrm{d}[2x-z]\\
&=-\underbrace{\mathrm{e}^{-z^2}\int_{-z}^{+z}y\mathrm{e}^{-y^2}\mathrm{d}y}_{=0}-z\mathrm{e}^{-z^2}
\int_{-z}^{+z}
\mathrm{e}^{-y^2}\mathrm{d}y\\
&=-\sqrt{\pi}z\mathrm{e}^{-z^2}\frac{1}{\sqrt{2\pi}}
\int_{-\sqrt{2}z}^{+\sqrt{2}z}\mathrm{e}^{-t^2/2}\mathrm{d}t\\
&=-\sqrt{\pi}z\mathrm{e}^{-z^2}\left[2\Phi\left(\sqrt{2}z\right)-1\right].
\end{aligned}
$$
Hence,
$$\int_0^z(1-\mathrm{e}^{-2(z-x)^2})4x\mathrm{e}^{-2x^2}\mathrm{d}x=1-\mathrm{e}^{-2z^2}-
\sqrt{\pi}z\mathrm{e}^{-z^2}\left[2\Phi\left(\sqrt{2}z\right)-1\right].$$

\end{document}